\documentclass{amsart}
\usepackage{amsfonts,amssymb,amstext}
\usepackage{url}
\usepackage{enumerate}
\usepackage{setspace}
\usepackage{tikz}

\newtheorem{theorem}{Theorem}[section]

\theoremstyle{definition}

\numberwithin{equation}{section}

\author{Acushla Sarswat}
\title[Nodal Line of a Second Eigenfunction of the Laplacian]{On the Nodal Line of a Second Eigenfunction of the Laplacian-Dirichlet in some Annular
Domains with Dihedral Symmetry}
\address{University of Mumbai, India}
\email{acushla.narayanan@gmail.com}
\thanks{The work of the author is 
supported by the Council of Scientific and Industrial Research, India via CSIR/JRF Award F. No. 09/019(0093)/2012-EMR-I.}


\begin{document}

\begin{abstract}
 Let $\Omega$ be a bounded annular $C^{1,1}$ domain in $\mathbb{R}^2$ which is left invariant
 under the action of the dihedral group $D_n$ of isometries of $\mathbb{R}^2$. We show that
 the nodal line of a second Dirichlet eigenfunction must intersect the boundary of $\Omega$,
 under suitable conditions on $\frac{\partial \;}{\partial \theta}$.  
\end{abstract}

\maketitle

\section{Introduction}

\onehalfspacing

In 1967, L. Payne \cite{Payne1967,Payne1973} conjectured that a second eigenfunction of the 
Laplacian with
Dirichlet boundary conditions cannot have a closed interior nodal curve for any 
bounded domain $\Omega \subset \mathbb{R}^2$, i.e. if $u_2$ is a solution of the problem
\cite[p.~6]{Henrot2006}
\begin{equation}\label{eigenvalue_problem}
 \left\{
 \begin{aligned}
  -\Delta u_2 &= \lambda_2 u_2 \hspace{0.3cm} \mbox{ in } \Omega, \\
          u_2 &= 0 \hspace{0.9cm} \mbox{ on } \partial \Omega
 \end{aligned}
 \right.
\end{equation}
where $\lambda_2$ is the second Dirichlet eigenvalue of $\Omega$, and the nodal line of $u_2$
is \[N = \overline{\{ x \in \Omega \; | \; u_2(x) = 0 \}},\] then we must have
\begin{equation}\label{statement_of_conjecture}
N \cap \partial \Omega \neq \varnothing. 
\end{equation}
L. Payne \cite{Payne1973} gave an explicit proof of this for domains
with smooth boundary which are
convex in $x$ and symmetric about the $y$-axis. In 1987, C.-S. Lin \cite{Lin1987} showed that  
it holds when $\Omega$ is symmetric under a rotation with angle
$2 \pi p/q$ where $p, q$ are positive integers. It has since been established 
\cite{Alessandrini1994,Jerison1991,Melas1992,Putter1990,YangG2013} that 
\eqref{statement_of_conjecture} holds true for all bounded, convex domains in $\mathbb{R}^2$ as well as for some simply-connected concave domains. 

However,
a characterization of all planar domains for which \eqref{statement_of_conjecture} holds is
now an open question, as counterexamples have been found \cite{HoffmannN1997} within the 
class of non simply-connected domains. 

In this paper we show that 
\eqref{statement_of_conjecture} holds within a class of annular,
dihedrally symmetric domains. The argument is an extension of L. Payne's proof \cite{Payne1973}.

\medskip

Let $(r,\theta)$ denote polar coordinates in $\mathbb{R}^2$
centered at the origin $(0,0)$ with
$r \geq 0, \theta \in ]-\pi,\pi]$. 
Our main result is the following:

\begin{theorem}\label{maintheorem}
Let $\Omega$ be an annular domain with $C^{1,1}$-boundary which is left invariant under the action of the dihedral group 
 $D_n$ of isometries of
 $\mathbb{R}^2$ generated by the rotation $\rho_n$ of $\mathbb{R}^2$ about
 some fixed point on the $x$-axis by
 an angle of $2\pi/n$, and the reflection $\sigma$ of $\mathbb{R}^2$
 about the $x$-axis. Assume that $\Omega$ is 
 contained in the region $-\frac{\pi}{2} < \theta < \frac{\pi}{2}$ and that the point $(0,0) \notin 
 \overline{\Omega}$.
Let $\partial \Omega = C_0 \mathaccent\cdot\cup C_1$ 
where the Jordan curves $C_0$ and $C_1$ are the outer 
and inner boundaries
of $\Omega$ respectively. 
Suppose further
that the following conditions hold:
\begin{enumerate}
\item $\frac{\partial}{\partial \theta}$ points outward from $\Omega$, i.e. 
$\left\langle \frac{\partial}{\partial \theta}(p), \nu(p) \right\rangle > 0$, 
$\forall p \in C_1 \cap \{ (r,\theta) \; | \; \theta < 0 \}$, (where $\nu(p)$ denotes the outward unit normal at $p$).
\item $\frac{\partial}{\partial \theta}$ points inward into $\Omega$, i.e. 
$\left\langle \frac{\partial}{\partial \theta}(p), \nu(p) \right\rangle < 0, \forall p \in C_0 \cap \{ (r,\theta) \; | \; \theta < 0 \}$.
\end{enumerate}
Then \eqref{statement_of_conjecture} holds, i.e. if $N=\overline{\{ x \in \Omega \; | \; u_2(x) = 0 \}}$ is the nodal line of a second eigenfunction $u_2$ of $\Omega$ then 
$N \cap \partial \Omega \neq \varnothing$.
\end{theorem}

Possible extensions of Theorem \ref{maintheorem} to a wider class of domains will be indicated in \S \ref{section2}. An example of a domain satisfying the above conditions is illustrated below.

\medskip

\begin{center}
\begin{tikzpicture}[scale=0.5]
\node (fig1) at (1,1) {Figure 1.};
  \path[draw, rounded corners=8pt] (-2.5,2.6698)--(-2.5,11.3301)--(5,7)--cycle;
\path[draw, rounded corners=8pt] (2,7)--(-1,8.7320)--(-1,5.2679)--cycle;
\end{tikzpicture} 
\end{center}

\section{Proof of Theorem \ref{maintheorem}}\label{section2}

If $u_2 \circ \sigma = -u_2$ 
there is nothing to 
show
because the nodal line of $u_2$ is then the intersection of $\Omega$
with the $x$-axis. 

Next consider the case where $u_2 \circ \sigma = u_2$.

\medskip

Suppose if possible that $N \subset \Omega$. 

By the Courant nodal domain theorem \cite{CourantH1953}, $\Omega \setminus N$ has two
connected components, $D_-$ and $D_+$, such that $u_2$ is strictly
negative in $D_-$ and
strictly positive in $D_+.$ 
Therefore as mentioned in \cite[\S2]{Payne1973} and as a consequence
of \cite[Theorem~2.5]{Cheng1976}, 
$N$ must be a loop.
Hence
the following cases arise.

\medskip
 
Case (i): $\partial D_- = N \mathaccent\cdot\cup C_1$
and $\partial D_+ = N \mathaccent\cdot\cup C_0.$

\medskip

Were $\frac{\partial u_2}{\partial \theta} \geq 0$ below the $x$-axis, it would imply $u_2 \geq 0$ below the $x$-axis and hence in all of $\Omega$, as $u_2=u_2 \circ \sigma$. 
Since this is not possible the set $R = \{ (r, \theta) \in \Omega \; \big|\; \theta < 0, \frac{\partial u_2}{\partial \theta}(r,\theta) < 0 \}$ must be
non-empty.

\medskip

Now $\Delta u_2 = -\lambda_2 u_2 \geq 0$ in $D_-$ and $u_2 = 0$ on $\partial D_-$. 
Since $\Omega$ has $C^{1,1}$-boundary it satisfies the interior sphere condition at
all points on $\partial \Omega$ by \cite[Theorem~1.0.9]{Barb2009}.
Hence by the Hopf lemma \cite{ProtterW1967}, $\frac{\partial u_2}{\partial \theta}(q) > 0, \forall q \in C_1 \cap \{ (r,\theta) \; | \; \theta < 0 \}$, given that $\frac{\partial}{\partial \theta}(q)$ 
is directed outwards from $D_-$ at these points. 

\medskip

Similarly $\Delta u_2 \leq 0$ in $D_+$ implies $\frac{\partial u_2}{\partial \theta}
(q) > 0$ for points $q \in C_0 \cap \{ (r,\theta) \; | \; \theta < 0 \}$. Hence
$\partial R \cap \partial \Omega$ is contained in the $x$-axis. Therefore $\frac{\partial u_2}{\partial \theta}=0$ on $\partial R$. 

\medskip

Since $\Omega$ has $C^{1,1}$-boundary, $\frac{\partial u_2}{\partial \theta} \in H^1(\Omega)$ \cite[Theorem~1.2.10]{Henrot2006},
and hence $\frac{\partial u_2}{\partial \theta} \in H^1(R \cup \sigma(R))$. As $\overline{R} \cap \overline{\sigma(R)}$ is contained in the $x$-axis,
$\frac{\partial u_2}{\partial \theta} \in H^1_0(R \cup \sigma(R))$.
Now $-\Delta(\frac{\partial u_2}{\partial \theta}) = \lambda_2 \frac{\partial u_2}{\partial \theta}$ in $\Omega$ in the weak sense.
Therefore
\[
 -\int_{R \cup \sigma(R)} \left(\frac{\partial u_2}{\partial \theta}\right)\Delta \varphi \mathop{dx} = \lambda_2 \int_{R \cup \sigma(R)} \left( \frac{\partial u_2}{\partial \theta} \right) \varphi 
 \mathop{dx}, \forall \varphi \in C_0^{\infty}(R \cup \sigma(R)).
\]
Since $C_0^{\infty}(R \cup \sigma(R))$ is dense in $H^1_0(R \cup \sigma(R))$ it follows
as a consequence of the Green's identity that
\[
 \lambda_2 = \frac{\int_{R \cup \sigma(R)} \lVert \nabla \left( \frac{\partial u_2}{\partial \theta}
 \right) \rVert^2 \mathop{dx}}{\int_{R \cup \sigma(R)} \left( \frac{\partial u_2}{\partial \theta}\right)^2 \mathop{dx}}
\]
Also, $\frac{\partial u_2}{\partial \theta}$ changes sign in $R \cup \sigma(R)$. Thus by the variational principle for Dirichlet eigenvalues $\lambda_2(R \cup \sigma(R))
< \lambda_2$ \cite[Formulae~1.35, Remark~1.2.4]{Henrot2006}.

\medskip

However this contradicts the domain monotonicity of Dirichlet eigenvalues \cite[\S1.3.2]{Henrot2006}. The argument is analogous when $\partial D_- = N \mathaccent\cdot\cup C_0$ and $\partial D_+ = N \mathaccent\cdot\cup C_1$.

\medskip

Case (ii): $\partial D_- = N$ and $\partial D_+ = \partial \Omega \mathaccent \cdot\cup N$. Let
\[w = \sum_{i=0}^{n-1} u_2 \circ \rho_n^i\]  and let
\[K = \mathop\cup_{i=0}^{n-1} \rho_n^i(\overline{D}_-).\] 
Then $\overline{D}_- \subset K$ and 
$\rho_n : K \to K$ is an isometry.
Therefore as $u_2$ is strictly positive in 
$\Omega \setminus K$, $w$ is also strictly positive in 
$\Omega \setminus K$. Since $w$ is non-zero, it is a second Dirichlet eigenfunction of $\Omega$. The nodal line of $w$ is contained in $K$ and therefore does not
intersect the boundary $\partial \Omega$. Now
\[
w \circ \sigma = \sum_{i=0}^{n-1} u_2 \circ \rho_n^i \circ \sigma = \sum_{i=0}^{n-1} u_2 \circ \sigma \circ \rho_n^{-i} = \sum_{i=0}^{n-1} u_2 \circ \rho_n^{-i} = w
\]
Also $w(q)=0 \implies w(\rho_n^i(q))=0, \forall i=0,\ldots,n-1$. Hence the nodal line of $w$ must be a loop encircling $C_1$. 

\medskip

Therefore the argument in case (i) applies to $w$ and this
leads to a contradiction. The situation is analogous when $\partial D_+ = N$ and $\partial D_- = \partial \Omega \mathaccent \cdot\cup N$. 

\medskip

Next we express $u_2$ as the sum 
$u_2 = v + w$ where $v = \dfrac{u_2 - u_2 \circ \sigma}{2}$ and 
$w = \dfrac{u_2 + u_2 \circ \sigma}{2}$. 
Suppose if possible that $N \subset \Omega$.

\medskip

Since $C_0$ is a compact subset of $\overline{\Omega} \setminus N$, there exists
an $\epsilon$-neighbourhood 
\[U_0 = \mathop\cup_{p \in C_0} B_{\epsilon}(p) \cap \overline{\Omega} \]
of $C_0$ in $\overline{\Omega} \setminus N$ for some $\epsilon > 0.$

\medskip

Since $U_0 \setminus C_0$ is connected, $u_2$ carries one sign in $U_0 \setminus C_0$, say ${u_2}_{|_{U_0 \setminus C_0}} > 0$. If $v \leq 0$ in $U_0 \cap \{ (r,\theta) \; | \; \theta \leq 0 \}$ then $w = u_2 - v \geq 0$ in $U_0 \cap \{ (r,\theta) \; | \; \theta \leq 0 \}$ and hence $w \geq 0$ in $U_0$
by symmetry. 
The same argument also works if $v \leq 0$ in $U_0 \cap \{ (r,\theta) \; | \; \theta \geq 0 \}$. 
As the nodal line of $v$ is the $x$-axis only these two possibilities exist.

\medskip

Thus if $q \in C_0$ is a point on the nodal curve of $w$ then $w$ does not change sign
in $B_{\epsilon}(q) \cap \Omega$, which is a contradiction. Hence
the nodal curve of $w$ does not intersect $C_0$. A similar argument shows that the nodal
curve of $w$ does not intersect $C_1$. But this contradicts what was shown earlier, because
$w=w \circ \sigma$ and $w$ is also a second eigenfunction of $\Omega$.

Therefore the nodal line of $u_2$ must intersect the boundary $\partial \Omega$. 
\hfill{$\Box$}

\bigskip

\noindent {\bf Concluding Remarks:}
\begin{enumerate}
\item The above proof will go through for a wider class of dihedrally symmetric annular
 domains provided we are able to show that $\frac{\partial u_2}{\partial \theta}$
 is in $H_0^1(R \cup \sigma(R))$. This holds for instance, when the boundary 
 of $\Omega$ is polygonal \cite[p.~IX]{Grisvard1992} and the rest of the
 conditions of Theorem \ref{maintheorem} are met. Conditions (1) and (2) of Theorem
 \ref{maintheorem} and the interior sphere condition are only required to hold almost everywhere on $\partial \Omega$ but we leave these considerations as a topic for future work as of now.
 \item The choice of the $x$-axis in Theorem \ref{maintheorem} is a matter of convenience. The proof goes through if polar coordinates are chosen in a way that the polar axis is an axis of reflection 
 of $\Omega$ and the pole is away from $\Omega$.
 \end{enumerate}

\medskip

\proof[Acknowledgements]
The author wishes to thank her thesis advisor Professor A. R. Aithal for his support as well as for discussing this problem and for pointing out a flaw in her initial argument, even
though this was not part of her thesis work. She also wishes to thank
Professor Jyotshana V. Prajapat of the University
of Mumbai for her support and for discussing the arguments presented in this paper.

\end{document}